\documentclass[12pt]{elsarticle}

\usepackage{tikz}
\usetikzlibrary{shapes}

\usepackage{fullpage}

\usepackage{amsmath,amsfonts,amssymb,amsthm,graphics,amscd}

\newcommand{\ceil}[1]{\left\lceil#1\right\rceil}
\newtheorem{theorem}{Theorem}[section]
\newtheorem{lemma}[theorem]{Lemma}

\newtheorem{conjecture}[theorem]{Conjecture}

\renewcommand\P{\mathbb{P}}
\def\ceil#1{\left\lceil #1 \right\rceil}
\def\Dln#1#2{D^{#1} \, \ln^{#2} D}
\def\De#1{D^{#1}}

\newcommand\Pry{Przyby{\l}o }

\begin{document}
\begin{frontmatter}

\title{Asymptotically optimal neighbor sum distinguishing total colorings of graphs}
\author[SL]{Sarah Loeb}
\ead{sloeb2@illinois.edu}
\author[JP]{Jakub Przyby{\l}o}
\ead{jakubprz@agh.edu.pl}

\author[YT]{Yunfang Tang\corref{cor1}}
\ead{tangyunfang8530@gmail.com}

\cortext[cor1]{corresponding author}
\address[SL]{Department of Mathematics, University of Illinois at Urbana--Champaign, 1409 W. Green St., Urbana, IL 61801, United States}
\address[JP]{AGH University of Science and Technology, al. A. Mickiewicza 30, 30-059 Krakow, Poland}
\address[YT]{Department of Mathematics, China Jiliang University, Xueyuan Road, Xiasha High Education District, Hangzhou, 310018, China}

\begin{abstract}
Given a proper total $k$-coloring $c:V(G)\cup E(G)\to\{1,2,\ldots,k\}$ of a graph $G$, we define the \emph{value} of a vertex $v$ to be $c(v) + \sum_{uv \in E(G)} c(uv)$. The smallest integer $k$ such that $G$ has a proper total $k$-coloring whose values form a proper coloring is the \emph{neighbor sum distinguishing total chromatic number} of $G$, $\chi''_{\Sigma}(G)$. Pil{\'s}niak and Wo{\'z}niak~(2013) conjectured that $\chi''_{\Sigma}(G)\leq \Delta(G)+3$ for any simple graph with maximum degree $\Delta(G)$. In this paper, we prove this bound to be asymptotically correct by showing that  $\chi''_{\Sigma}(G)\leq \Delta(G)(1+o(1))$. The main idea of our argument relies on Przyby{\l}o's proof (2014) regarding neighbor sum distinguishing edge-colorings.\\

\end{abstract}

\end{frontmatter}

\section{Introduction}

A \emph{proper $k$-coloring} of a graph $G$ is a function $c: V(G) \to [k]$ such that $c(u) \neq c(v)$ whenever $uv \in E(G)$, where $[k]$ denotes $\{1,\ldots,k\}$. The least $k$ such that $G$ has a proper $k$-coloring is the \emph{chromatic number} of $G$, $\chi(G)$. Clearly $\chi(G) \le \Delta(G) +1$, where $\Delta(G)$ is the maximum degree of $G$. A \emph{proper $k$-edge-coloring} of $G$ is a function $c:E(G) \to [k]$ such that $c(uv) \neq c(uw)$ when $v \neq w$. The least $k$ such that $G$ has a proper $k$-edge-coloring is the \emph{edge-chromatic number} of $G$, $\chi'(G)$. For a $k$-edge-coloring $c$, define the \emph{value} $s_c(v)$ of a vertex $v$ by $\sum_{u \in N(v)} c(uv)$. A $k$-edge-coloring $c$ of $G$ is a \emph{proper edge-weighting} if $s_c$ forms a proper coloring of $G$. The least $k$ such that $G$ has a proper $k$-edge-coloring that is a proper edge-weighting is the \emph{neighbor sum distinguishing edge-chromatic number} of a graph, denoted $\chi'_{\Sigma}(G)$. Note that this graph parameter is well defined only for graphs with no isolated edges.

The study of this type of graph invariants, focused on distinguishing vertices by their associated sums of colors of e.g. their incident edges was initiated in~\cite{Chartrand},
where the parameter called the \emph{irregularity strength} of a graph was introduced. Within this, motivated by the trivial fact that no irregular graph exists (with 1-vertex exception) and with reference to research in~\cite{ChartrandErdosOellermann} of Chartrand, Erd\H{o}s and Oellermann on possible alternative definitions of such a graph, Chartrand et al. tried to propose a measure of irregularity of a given graph $G$ by means of multiplying edges of $G$ in order to produce an irregular multigraph of $G$, see~\cite{Chartrand} for details. The same construction was also expressed in terms of colorings, where integer colors of edges of a graph corresponded to multiplicities of the edges in a corresponding multigraph, while the sums of these colors from the edges incident with a given vertex -- to the degree of this vertex in the multigraph. The irregularity strength is a well studied graph invariant, but it also gave rise to a whole discipline, which might be referred to as additive graph labellings or more generally -- vertex distinguishing graph colorings, see e.g.~\cite{Gallian_survey}, including extensive study of these.

For the edge-coloring parameters, $\chi'_{\Sigma}(G) \ge \chi'(G) \ge \Delta(G)$, while by Vizing's Theorem, $\chi'(G) \le \Delta(G)+1$ for every graph $G$. Flandrin, Marczyk, Przyby{\l}o, Sacl\'{e}, and Wo{\'z}niak~\cite{FMPSW} conjectured that:

\begin{conjecture}[\cite{FMPSW}] \label{conj:edge}
If $G$ is a connected graph with at least three vertices other than $C_5$, then $\chi'_{\Sigma}(G) \le \Delta(G) + 2$.
\end{conjecture}

Przyby{\l}o~\cite{P} proved an asymptotically optimal upper bound for graphs with large maximum degree. Specifically, he showed:

\begin{theorem}[\cite{P}] \label{thm:P}
If $G$ is a connected graph with $\Delta(G)$ sufficiently large, then $\chi'_{\Sigma}(G) \le \Delta(G) + 50 \Delta(G)^{5/6} \ln^{1/6} \Delta(G)$.
\end{theorem}

A \emph{proper total $k$-coloring} of $G$ is a function $c: V(G) \cup  E(G) \to [k]$ such that $c$ restricted to $V(G)$ is a proper coloring, $c$ restricted to $E(G)$ is a proper edge-coloring, and such that the color on each vertex is different from the color on its incident edges. The least number of colors in such a coloring of $G$ is denoted by $\chi''(G)$. For a proper total $k$-coloring $c$, define the \emph{value} $s_c(v)$ of a vertex $v$ by $c(v) + \sum_{u \in N(v)} c(uv)$. A proper total $k$-coloring $c$ of $G$ is a \emph{proper total weighting} if $s_c$ is a proper coloring of $G$. The least $k$ such that $G$ has a proper total $k$-coloring that is a proper total weighting is the \emph{neighbor sum distinguishing total chromatic number}
of $G$, denoted $\chi''_{\Sigma}(G)$. Clearly, $\chi''_{\Sigma}(G) \ge \chi''(G) \ge \Delta(G)+1$.
On the other hand, the famous Total Coloring Conjecture, that has eluded mathematicians for half a century presumes that $\chi''(G)\leq\Delta(G)+2$ for every graph $G$. This was independently posed by Vizing~\cite{Vizing2} and Behzad~\cite{Behzad}.
Thus far, it has been confirmed up to a (large) additive constant by means of the probabilistic method, see~\cite{MolloyReedTotal}. Despite that, Pil{\'s}niak and Wo{\'z}niak~\cite{PW} daringly conjectured the following.

\begin{conjecture}[\cite{PW}] \label{conj:total}
If $G$ is a 
graph with maximum degree $\Delta(G)$, then $\chi''_{\Sigma}(G)\leq \Delta(G)+3$.
\end{conjecture}

Pil{\'s}niak and Wo{\'z}niak~\cite{PW} proved that Conjecture~\ref{conj:total} holds for complete graphs, cycles, bipartite graphs and subcubic graphs. Using the Combinatorial Nullstellensatz, Wang, Ma, and Han~\cite{WMH} proved that the conjecture holds for triangle-free planar graphs with maximum degree at least $7$.  Dong and Wang~\cite{DW} showed that Conjecture~\ref{conj:total} holds for sparse graphs, and Li, Liu, and Wang~\cite{LLW} proved that the conjecture holds for $K_4$-minor free graphs. Li, Ding, Liu, and Wang~\cite{LDLW} also confirmed Conjecture~\ref{conj:total} for planar graphs with maximum degree at least 13. Finally, Xu, Wu, and Xu~\cite{XWX} proved $\chi''_{\Sigma}(G)\leq \Delta(G)+2$ for graphs $G$ with $\Delta(G) \ge 14$ that can be embedded in a surface of nonnegative Euler characteristic.

By modifying Przyby{\l}o's proof that Conjecture~\ref{conj:edge} is asymptotically correct for graphs with large maximum degree, we confirm in this paper that Conjecture~\ref{conj:total} is also asymptotically correct by showing the following. 

\begin{theorem}\label{thm:total} If $G$ is a 
graph with $\Delta(G)$ sufficiently large, then $\chi''_{\Sigma}(G) \le \Delta(G) + 50 \Delta(G)^{5/6} \ln^{1/6} \Delta(G)$.
\end{theorem}

\section{Ideas}

We color the vertices of the graph first and produce an edge-coloring such that the combined total coloring is a proper total weighting. For a coloring $g$ and an edge-coloring $h$, let $g \sqcup  h$ be the total coloring produced by combining $g$ and $h$.

Our main work is in producing the desired edge-coloring. Our Lemma~\ref{lem:long} serves a similar purpose to Lemma~6 of Przyby{\l}o~\cite{P}. Lemma~\ref{lem:long} guarantees (not necessarily proper) colorings $c_1$ and $c_2$ of the vertices and edges respectively such that the colors are roughly evenly distributed. These colorings are used to produce an initial (also improper) edge-coloring $c'$ by setting $c'(uv) = c_1(u) + c_1(v) + c_2(uv)$. Statement~(3) of Lemma~\ref{lem:long} guarantees that the colors used by $c'$ are also roughly evenly distributed. Finally, statement~(4) of Lemma~\ref{lem:long} will be used to guarantee that the final values for the vertices form a proper coloring. The proof uses the Lov\'{a}sz Local Lemma and the Chernoff Bound in the forms below.

\begin{theorem}[Lov\'{a}sz Local Lemma~\cite{LLL}] \label{thm:LLL}
 Let $A_1,\ldots,A_n$ be events in an arbitrary probability space. Suppose that each event $A_i$ is mutually independent of a set of all but at most $D$ others of these events, and that $\P(A_i) \le p$ for all $1 \le i \le n$. If $ep(D+1) \le 1$, then $\P(\bigcap_{i=1}^n \overline{A_i}) > 0$.
\end{theorem}

\begin{theorem}[Chernoff Bound~\cite{Cher}] \label{thm:CB}
If $0 \le t \le np$, then
\[ \P(|\operatorname{BIN}(n,p) - np| > t) < 2e^{-t^2/3np}, \]
where $\operatorname{BIN}(n,p)$ is a binomial random variable with $n$ independent trials having success probability $p$.
\end{theorem}

For a vertex $v$, let $d(v)$ be the degree of $v$.

\begin{lemma} \label{lem:long}
Let $H$ be a graph of maximum degree at most $D$. For any coloring $c_1$ of $V(H)$ and a vertex $v \in V(H)$ with $d=d(v)$, let
\[S(v) =   \left( \ceil{D^{2/3} \ln^{1/3} D} + \ceil{D^{1/2}} \right) d c_1(v) + R(d,D), \]
where $R$ is a function of two variables.

If $D$ is sufficiently large, then there exist colorings $c_1: V(H) \to \left[ \ceil{D^{1/6}\ln^{-1/6} D} \right]$ and $c_2: E(H) \to \left[ \ceil{D^{1/3}\ln^{-1/3} D} \right]$ such that if $S(v) \le D^2$ for all $v \in V(H)$, then for every vertex $v$: \begin{enumerate}
\renewcommand{\theenumi}{(\arabic{enumi})}
\item if $d \ge \Dln{5/6}{1/6}$, then the number of neighbors $u$ of $v$ having any given color $c_1(u)$ is within $3 \Dln{5/12}{7/12}$ of $d\left(\ceil{\Dln{1/6}{-1/6}}\right)^{-1}$;
\item if $d\ge \Dln{5/6}{1/6}$, then the number of edges $e$ incident with $v$ having any given color $c_2(e)$ is within $3 \Dln{1/3}{2/3}$ of $d \left(\ceil{\Dln{1/3}{-1/3}}\right)^{-1}$;
\item for every integer $c$ with $3 \le c \le 2 \ceil{\Dln{1/6}{-1/6}} + \ceil{\Dln{1/3}{-1/3}}$, the number of edges $uv$ incident with $v$ such that $c_1(u) + c_1(v) + c_2(uv) = c$ is at most $2 \Dln{1/3}{1/3}$ if $d< \De{2/3}$, or at most $\Dln{2/3}{1/3} + 3 \Dln{1/3}{2/3}$ otherwise;
\item if $d\ge 3 \Dln{5/6}{1/6}$, then for every integer $\alpha > 0$, the number of neighbors $u$ of $v$ such that $\frac{d}{2} \le d(u) \le 2 d$ and $S(u) \in I_{d,\alpha,D}$ is at most $d\left(\Dln{1/6}{-1/6}\right)^{-1} + 3\Dln{5/12}{7/12}$, where
\[ I_{d,\alpha,D} = \left( (\alpha-1) \frac{d}{2} D^{2/3} \ln^{1/3} D, \alpha \frac{d}{2} D^{2/3} \ln^{1/3} D \right]. \]
\end{enumerate}
\end{lemma}

{
\renewenvironment{proof}{\emph{Proof (Sketch).}}{\qed}
\begin{proof}
The only difference between Lemma~6 of~\cite{P} and Lemma~\ref{lem:long} above is that \Pry writes:
\[ S(v) = \left( \ceil{\Dln{2/3}{1/3}} + 4 \ceil{\Dln{1/3}{2/3}} \right) d c_1(v) + R(d,D)\]
and we have
\[S(v) =   \left( \ceil{\Dln{2/3}{1/3}} + \ceil{D^{1/2}} \right) d c_1(v)+ R(d,D).\]
 Our change to the lower order term in $S(v)$ comes from an increase in the number of possible colors to be used on the edges and is made to accommodate a total coloring. However, the proof is unchanged, so we omit the details and give only a sketch of the ideas.
 
Start with colorings $c_1$ and $c_2$ where the color assigned to each vertex and edge is chosen independently and uniformly at random. For each vertex $v$, define four events corresponding to $v$ violating each of (1), (2), (3), and (4). For each bad event $E$, the Chernoff bound shows that the probability of $E$ is less than $D^{-5/2}$. Since all events for a vertex $v$ are mutually independent of those corresponding to vertices at distance at least three from $v$, each event is mutually independent of all but at most $3+4D^2$ other events. Finally, since
\[ e D^{-5/2} (4+4 D^2) < 1, \]
the Lov\'{a}sz Local Lemma implies that there is some selection of $c_1$ and $c_2$ such that none of the bad events occurs.
\end{proof} }

To form a total coloring, we start with a coloring of the vertices and extend it to a total coloring. To guarantee that the total coloring is proper, we use a result of Molloy and Reed~\cite{MR}. A \emph{list assignment} $L$ for $E(G)$ assigns to each edge $e$ a list $L(e)$ of permissible colors. Given a list assignment $L$ for the edges of $G$, if a proper edge-coloring $c$ can be chosen so that $c(e)\in L(e)$ for all $e\in E(G)$, then we say that $G$ is \emph{$L$-edge-colorable}. The \emph{list edge-chromatic number} $\chi_{\ell}'(G)$ of $G$ is the least $k$ such that $G$ is $L$-edge-colorable for any list assignment $L$ satisfying $|L(e)|\ge k$ for all $e\in E(G)$.

\begin{theorem}[Molloy and Reed~\cite{MR}] \label{thm:listedge}
There is a constant $k$ such that $\chi_{\ell}'(G) \le \Delta(G) +k \Delta(G)^{1/2}(\log \Delta(G))^4$ for every graph $G$.

\end{theorem}

\section{Proof of Theorem~\ref{thm:total}}

We first give an outline of the proof.

Suppose $G$ is a graph and $g:V(G) \to [\Delta(G)+1]$ is a proper coloring of $G$. We will produce a proper edge-coloring $h$ such that $g \sqcup  h$ is a proper total weighting. Let $M$ be a maximal matching in $G$. Producing $h$ takes three steps: the first two steps focus on producing an edge-coloring of $G-M$, and the final step assigns colors to $M$.

More specifically, in Step 1, we use Theorem~\ref{thm:listedge} to define an edge-coloring $h_1$ for $E(G)-M$ such that $g \sqcup  h_1$ is a proper total coloring of $G-M$. In Step 2, we modify $h_1$ to obtain an edge-coloring $h_2$ on $E(G)-M$ so that $h_2 \sqcup  g$ is a proper total coloring and $s_{g \sqcup  h_2}(u) \neq s_{g \sqcup  h_2}(v)$ whenever $uv \in M$. In Step 3, we extend $h_2$ to $M$ to obtain a coloring $h$ of $E(G)$ such that $g \sqcup  h$ is a proper total coloring that is also a proper total weighting of $G$.

\begin{proof}
Let $G$ be a graph with maximum degree $D$. Let $M$ be a maximal matching in $G$, and define $G'$ by setting $V(G') = V(G)$ and $E(G') = E(G)-M$. Let $g:V(G) \to [D+1]$ be a proper coloring of $G$ (and thus of $G'$).

Let $C_V = \ceil{D^{1/6} \ln^{-1/6} D}$ and  $C_E = \ceil{D^{1/3} \ln^{-1/3} D}$. These are the numbers of colors used in the coloring $c_1$ and the edge-coloring $c_2$ guaranteed by Lemma~\ref{lem:long}. Let $d_0 = D^{5/6} \ln^{1/6} D$ so that $3d_0$ is the the degree threshold in Lemma~\ref{lem:long}~(4).  Let $C_M = \ceil{47 d_0}$; we will color $M$ from $[C_M]$. The dominant term in the ``stretch factor'' used to produce a proper edge-coloring is $\ceil{D^{2/3} \ln^{1/3} D}$, which we abbreviate as $C$.

\medskip

\textbf{Step 1:} The coloring $h_1$ is defined in several phases. Our argument follows that of Sections~5.1 and~5.2 in~\cite{P}, with modifications to produce a total coloring rather than an edge-coloring.

Let $c_1:V(G') \to [C_V]$ and $c_2:E(G') \to [C_E]$ be the colorings of $V(G')$ and $E(G')$ guaranteed by Lemma~\ref{lem:long}, where the function $R(d,D)$ will be specified later.

Assign $uv \in E(G')$ a tentative color $c'(uv)$ which we define by
\[ c'(uv) = \left[ c_1(u) + c_1(v) + c_2(uv) \right] \left( C + \ceil{\De{1/2}} \right) + C_M. \]
This coloring is not a proper edge-coloring. However, by Lemma~\ref{lem:long}~(3), the colors are distributed so that we will be able to modify them to produce a proper edge-coloring $h_1$.  Note that
\begin{equation}\label{c_prime_ineq}
\ceil{47d_0}<c'(uv)\le (2C_V + C_E)(C + \ceil{\De{1/2}}) + C_M\le D+49d_0+o(d_0).
\end{equation}
The colors $1$ through $C_M$ are not used until Step 3, when they are used on $M$.

For each $\beta \in \{3,\ldots,2 C_V + C_E\}$, the set of integers from $\beta \left( C + \ceil{\De{1/2}} \right) + C_M$ to $(\beta + 1) \left( C + \ceil{\De{1/2}} \right) + C_M -1$ will be called a \emph{$\beta$-palette}. Note that at this point of the construction, only the smallest member of each palette may appear as a color of an edge of $G'$. Each edge $uv$ in $E(G')$
is now assigned a $\beta$-palette with $\beta=c_1(u)+c_1(v)+c_2(uv)$. Note that $c'(uv)$ belongs to such $\beta$-palette. This will also hold throughout the construction for
$h_1(uv)$ and $h_2(uv)$. We will now define $h_1(e) = c'(e) + a_1(e)$, where $a_1(e)$ specifies which element from the palette associated with $e$ is assigned to $e$. To this end, let $P = \{0,\ldots,C+\ceil{D^{1/2}}-1\}$. We divide $P$ into a lower and upper portion with $P^- = \{0,\ldots,C+\ceil{D^{1/2}/2}\}$ and $P^+ = \{C+\ceil{D^{1/2}/2}+1,\ldots,C+\ceil{D^{1/2}}-1\}$.

In our specification of the final coloring $h$, we will have $h(e) - c'(e) \le C + \ceil{\De{1/2}}$, so if $D$ is sufficiently large, then by~(\ref{c_prime_ineq}),
\[ h(e) \le D+49d_0+o(d_0) + C + \ceil{\De{1/2}}< D + 50 \Dln{5/6}{1/6} = D + o(D). \]

To choose $a_1(e)$, we first specify a list assignment and then use Theorem~\ref{thm:listedge}. For every consecutive $\beta\in\{3,\ldots,2C_V+C_E\}$ we now proceed as follows. Let $G_\beta$ be the spanning subgraph of $G'$ with $E(G_\beta) = \{e \in E(G'):c'(e) = \beta(C+\lceil D^{1/2}\rceil)+C_M\}$. By Lemma~\ref{lem:long}~(3), $\Delta(G_\beta) \le \Dln{2/3}{1/3} + 3 \Dln{1/3}{2/3}$. For an edge $uv\in E(G_\beta)$, let 
$T_\beta(uv) = \left\{g(u) - c'(uv), g(v) - c'(uv)\right\}$. To guarantee that edges receive colors distinct from the colors of their endpoints, let
$L_\beta(uv) = P^- - T_\beta(uv)$. For $D$ sufficiently large, we thus have
\[ |L_{\beta}(uv)| \ge C + \ceil{\De{1/2}/2} -2 \ge \Delta(G_\beta) + k \Delta(G_\beta)^{1/2} \log^4(\Delta(G_\beta)), \]
where $k$ is the constant from Theorem~\ref{thm:listedge}.
Let $a_{1,\beta}$ be an $L_\beta$-edge-coloring for $G_\beta$ guaranteed by Theorem~\ref{thm:listedge}.
For every $uv\in E(G_\beta)$ we then set $a_1(uv)=a_{1,\beta}(uv)$.

Note that in $h_1$, only the lower portion of the elements from each palette is used. The remaining colors in the palettes are used in Step 2.

The definition of $a_1$ guarantees that under $h_1$ no color is used on two incident edges. Thus $h_1$ is a proper edge-coloring. Furthermore,
for every $uv\in E(G')$ we have $h_1(uv) \notin \{g(u),g(v)\}$, so $g \sqcup  h_1$ is a proper total coloring of $G'$.

\medskip
\noindent \noindent \textbf{Step 2:} This step has two phases, with no substantial difference between our argument and that of Section~5.3 in~\cite{P}.

For every vertex of degree at least $\De{2/3},$ we randomly select a neighboring edge from $G'$. Using the Lov\'{a}sz Local Lemma, we can choose these edges so that $d_H(v) -1 \le 2 \De{1/3} \le \frac{1}{2}(\frac{1}{2} \De{1/2} - 6)$ for every $v \in V(G)$, where $H$ is the subgraph induced by the chosen edges. We examine the edges of $H$ one by one (in any order). When we reach the last edge $uv$ of $H$ incident with any edge (or two edges) of $M$, we modify the color on $uv$ if necessary so that the ends of its incident edge (or two edges) from $M$ have distinct values assigned afterwards. In order to achieve the described goal, we pick $a_2'(uv) \in P^+ - T_\beta(uv)$ (where $uv$ is assigned the $\beta$-palette) so that replacing $a_1(uv)$ with $a_2'(uv)$ for such an edge $uv$ preserves properness of the total coloring as well. The above bound on the degrees in $H$ makes this possible. Let $a_2'(e) = a_1(e)$ for all other (unmodified) edges of $G'$, and let $h_2'(e) = c'(e) + a_2'(e)$ for every $e\in E(G')$.

Let $M'$ be the set of edges in $M$ whose endpoints currently do not have distinct colors. For each edge $uu' \in M'$, we pick an edge $uw \in E(G')$ incident with $uu'$ in $G$. By the choice of $h_2'$, we have $d(u) < \De{2/3}$ (as otherwise $u$ and $u'$ would be sum distinguished due to the construction of $h'_2$ from the previous paragraph), so by Lemma~\ref{lem:long}~(3), there are at most $2 D^{1/3} \ln^{1/3} D$ members of the palette assigned to $c'(uw)$ on the edges incident with $u$ and at most $C + 3\Dln{1/3}{2/3}$ members of this palette on the edges incident with $w$. We may thus easily pick $a_2(uw) \in  P - T_\beta(uw)$ so that (after replacing $a'_2(uw)$ with $a_2(uw)$) the values of $u$ and $u'$ are different and that if $w w' \in M$,  then the values of $w$ and $w'$ are also different (and so that the total coloring remains proper).

For all other edges $e$ of $G'$, set $a_2(e) = a_2'(e)$ and let $h_2(e) = c'(e) + a_2(e)$.

\medskip
\noindent \noindent \textbf{Step 3:} This step follows the argument from Section~5.4 from~\cite{P}. We present the argument with the computations omitted and simply emphasize the change made to accommodate a total coloring.

Before we define $h$, we need to know (roughly) the current value of the vertices. For any $v \in V(G)$ we have
\[ s(v) = g(v) + \sum_{u \in N(v)} h_2(uv). \]
Let $v$ be a vertex with $d=d(v) \ge d_0$. Using Lemma~\ref{lem:long}~(1) and~(2), we may break $s(v)$ into a dominant term $S(v)$ and an error term $F(v)$. Specifically, we write $s(v) = S(v) + F(v)$ where
\begin{equation}\label{SofVeq}
S(v) = g(v) + d C_M + \left( C + \ceil{\De{1/2}} \right) \\
\times \left[ d c_1(v) + \frac{d}{C_V} \binom{C_V +1}{2} + \frac{d}{C_E} \binom{C_E+1}{2} \right]
\end{equation}
and
\begin{equation}\label{FofVeq}
F(v) = \sum_{u \in N(v)} a_2(uv) + \left(C + \ceil{\De{1/2}} \right) \times \left[ f_1(v) \binom{C_V + 1}{2} + f_2(v) \binom{C_E + 1}{2} \right],
\end{equation}
where $|f_1(v)|\le 3 \Dln{5/12}{7/12}$ and $|f_2(v)|\le 3 \Dln{1/3}{2/3}$ result from the error terms in the distribution of colors from Lemma~\ref{lem:long}~(1) and~(2).

We define $R$ so that this $S(v)$ is the one needed to apply Lemma~\ref{lem:long}. Note that $S(v) \le \De{2}$ when $D$ is sufficiently large. Additionally, notice that by~(\ref{SofVeq}) and~(\ref{FofVeq}), $S(v)\sim\frac{dD}{2}$, while $|F(v)| \le \left( \frac{5}{2} + o(1) \right) \Dln{5/3}{1/3}$. Thus \begin{equation}\label{3858ineq}
\frac{3}{4} \frac{d D}{2} < s_{g \sqcup  h_2}(v) < \frac{5}{4} \frac{d D}{2}
\end{equation}
when $D$ is sufficiently large (and this will not change if we later increase the value of every vertex by an irrelevant additive non-negative factor of at most $C_M$ while choosing colors for the edges in $M$).

Each vertex $v$ with $d < 3 d_0$ automatically has fewer than $\frac{1}{2}C_M-1$ neighbors whose value can equal the value of $v$. Consider a vertex $v$ with $d \ge 3 d_0$. If $u$ is a neighbor of $v$ satisfying $d(u) \ge d_0$ and $d(u) \notin [\frac{d}{2},2d]$, then
$s_{g \sqcup  h_2}(u)<\frac{5}{8}d(u)D<\frac{5}{16}dD$ (if $d(u)<\frac{d}{2}$) or $s_{g \sqcup  h_2}(u)>\frac{3}{8}d(u)D>\frac{3}{4}dD$ (if $d(u)>2d$),
hence the values of $u$ and $v$ differ by more than $C_M$ under $g \sqcup  h_2$ due to inequalities~(\ref{3858ineq}). Additionally, by~(\ref{c_prime_ineq}), every neighbor $u$ of $v$ with $d(u)<d_0$ will have value at most $d(u)(D+50d_0)+(D+1)<\frac{9}{8}d_0D\leq \frac{3}{8}dD<s_{g \sqcup  h_2}(v)$ by~(\ref{3858ineq}). Thus the only neighbors of $v$ whose values could eventually equal that of $v$ have degrees in $[\frac{d}{2},2d]$. The upper bound for $|F(v)|$ thus shows that any neighbor $u$ of $v$ that eventually could have value equal to the value of $v$ has $S(u)$ in one of the intervals $I_{d,\alpha,D}$ for at most $2[(10+o(1))\frac{D}{d}+1]+1$ consecutive values of $\alpha$. Thus, Lemma~\ref{lem:long} guarantees that there are at most $(23 + o(1)) d_0$ neighbors of $v$ whose values could eventually equal to that of $v$.

We may thus extend $h_2$ to $h$ using colors in $[C_M]$ on the edges in $M$ so that $g \sqcup  h$ is a proper total coloring that is a proper total weighting of $G$.

\end{proof}

\section{Acknowledgments}
The authors would like to thank Douglas B. West for guidance and for editing suggestions.
The second author was supported by the National Science Centre, Poland, grant no. 2014/13/B/ST1/01855 and partly
supported by the Polish Ministry of Science and Higher Education.
 The third author performed research at University of
Illinois and was supported by Chinese Scholarship Council.

\bibliography{NSDTotalColoring}

\end{document}